\documentclass[12pt, a4paper]{article}

\usepackage{amsmath}\multlinegap=0pt
\usepackage[latin1]{inputenc}
\usepackage{amsthm}
\usepackage{amssymb}
\usepackage[francais,english]{babel}
\numberwithin{equation}{section}
\theoremstyle{plain}
\newtheorem{thm}{Theorem}[section]

\newtheorem{prop}[thm]{Proposition}
\newtheorem{lem}[thm]{Lemma}
\newtheorem{defi}[thm]{Definition}

\theoremstyle{definition}

\theoremstyle{remark}
\newtheorem{rem}[thm]{Remark}
\newtheorem{ex}[thm]{Exemple}

\newcommand{\ve}{\varepsilon}

\newcommand{\lcal}{\mathcal{L}}

\newcommand{\R}{\mathbb{R}}

\newcommand\huno{{\mathcal H}^1}

\newcommand\N{{\mathbb N}}

\newcommand\pref[1]{(\ref{#1})}

\let \eps\varepsilon

\newcommand{\res}{\mathop{\hbox{\vrule height 7pt width .5pt depth 0pt \vrule height .5pt width 6pt depth 0pt}}\nolimits}

\newcommand\A{{\cal A}}

\newcommand{\weakto}{\rightharpoonup}

\newcommand\M{{\cal M}}
\newcommand\HH{{\cal H}}

\newcommand\clos{\overline}

\newcommand\omb{\clos{\Omega}}

\def\<#1,#2>{\left<#1,#2\right>}

\newcommand\cb{\overline{c}}
\newcommand\cbxi{\cb_{\xi}}
\newcommand\tilC{\widetilde{C}}
\newcommand\tilQ{\widetilde{Q}}
\newcommand\tilsig{\widetilde{\sigma}}
\newcommand\QQ{{\cal Q}(\mu_0,\mu_1)}
\newcommand\QQq{{\cal Q}^q(\mu_0,\mu_1)}
\newcommand\QQqg{{\cal Q}^q(\ogm)}
\newcommand\Leb{{\cal L}}
\newcommand\ogm{{\overline\gamma}}
\newcommand\op{{\overline p}^{x,y}}
\newcommand\ovp{{\overline p}}
\newcommand\oxi{{\overline\xi}}
\newcommand\oQ{{\overline Q}}
\newcommand\oc{{\overline c}}
\newcommand\e{\varepsilon}
\newcommand\ws{\stackrel{*}{\rightharpoonup}}

\title {Optimal transportation with
traffic congestion and Wardrop equilibria}
\author {G.~Carlier, C.~Jimenez \thanks{\scriptsize\ CEREMADE, UMR CNRS 7534, Universit\'e Paris IX
Dauphine, Pl. de Lattre de Tassigny, 75775 Paris Cedex 16, FRANCE
\texttt{carlier@ceremade.dauphine.fr, jimenez@ceremade.dauphine.fr }.},  F.~Santambrogio
\thanks{\scriptsize\ Scuola Normale Superiore, Classe di Scienze, Piazza dei Cavalieri 7, 56126, Pisa, ITALY \texttt{f.santambrogio@sns.it}. }}

\begin{document}

\maketitle

\begin{abstract}
In the classical Monge-Kantorovich problem, the transportation cost only depends on the amount of mass sent from sources to destinations and not on the paths followed by this mass. Thus, it does not allow for congestion effects. Using the notion of traffic intensity, we propose a variant taking into account congestion.  This leads to an optimization problem posed on a set of probability measures on a suitable paths space. We establish existence of minimizers and give a characterization. As an application, we obtain existence and variational characterization of equilibria of Wardrop type in a continuous space setting.
\end{abstract}


\clearpage

\section{Introduction}

Given two mass distributions $\mu_0$ and $\mu_1$ on $\R^d$ with equal total mass, the classical Monge-Kantorovich problem consists in finding  transportation plans (i.e. measures on $\R^d\times \R^d$ having $\mu_0$ and $\mu_1$ as marginals) with minimal average transportation cost. This old  problem  can be traced back to Monge \cite{monge}. It has received a lot of attention in the recent years since the pathbreaking paper of Brenier \cite{bre} who solved the case of a quadratic transportation cost.  We refer to the book of Villani \cite{villani}, the lecture notes of Ambrosio \cite{ambrosio} and the references therein for a recent account of this rich mathematical theory and its numerous applications.  

\smallskip

An interesting  case is when the transportation cost is given by a conformally flat Riemannian distance:
\[d_g(x,y):=\inf \left\{\int_0^1 g(\sigma(t))\vert \dot{\sigma}(t)\vert \; dt  \; :\; \sigma(0)=x, \sigma(1)=y \right\}\]
the corresponding Monge-Kantorovich problem then reads as
\begin{equation}\label{mkg}
\inf_{\gamma\in\Pi(\mu_0,\mu_1)} \int_{\R^d\times \R^d} d_g(x,y)d\gamma(x,y)\end{equation}
where $\Pi(\mu_0,\mu_1)$ is the set of transportation plans. However, from a traffic planning point of view,  problem \pref{mkg}  is not totally realistic. On the one hand,  \pref{mkg} is \emph{path-independent}: the total transportation cost  only depends  on the amount of mass transported from the sources $x$ to the destinations $y$ and not on the paths followed by this mass. Put differently, in \pref{mkg}, individual's travelling strategies are irrelevant. On the other hand, \pref{mkg} does not take into account congestion effects i.e. the possibility that the cost $g(x)$ for passing through the point $x$ depends on ``how crowded" small neighbourhoods of $x$ are. This idea can be made precise thanks to the notion of \emph{traffic intensity} associated to a probability measure on a suitable set of paths.  This notion of \emph{traffic intensity} (see paragraph \ref{tcm})  is the path-dependent analogue  of the well-known notion of \emph{transport density} in Monge's problem (see Bouchitt\'e, Buttazzo and Seppecher \cite{bou}, Bouchitt\'e and Buttazzo \cite{bb}, Caffarelli, Feldman and McCann \cite{cfm}). With this notion at hand, we propose an optimal transportation problem with congestion. This variant of \pref{mkg}  takes the form of a relatively simple convex optimization problem posed on a set of probability measures on a suitable path space. We obtain existence  of minimizers (theorem \ref{existence}) and a characterization (Theorem \ref{caractopt}). 
\smallskip

Researchers in the field of applied traffic modelling have long emphasized the role of congestion in networks. In the early 50's, Wardrop (see \cite{wardrop}) considered the situation where a large number of vehicles have to go from one location to another, connected by a finite number of different roads. Each vehicle has to choose one road (or a probability on the set of possible roads) to minimize some transportation cost which depends not only on the road chosen but also on the total flow of vehicles on this road. Wardrop gave a minimal stability requirement for transportation strategies: the  cost of every actually used road should be equal or less than that which would be experienced by a single vehicle on any unused road. This natural equilibrium concept has been very popular since its introduction because of applications to networks of course but also due to the development of non-cooperative game theory in the meanwhile. From the best of our knowledge, the study of Wardrop equilibria have mainly been restricted to the case where admissible roads are given by a finite graph. A secondary contribution of the present paper is to introduce an equilibrium concept of Wardrop type in a continuous state setting, to prove the existence of such equilibria and to relate it to the optimal transportation problem with congestion (theorem \ref{cwe}).

\section{Optimal transportation with congestion}

\subsection{Notations}

Given a locally compact separable metric space $X$, we will denote respectively by  $\M_+(X)$ and $\M_+^1(X)$ the set of positive and finite Radon measures on $X$ and the  set of Radon probability measures on $X$. If $X$ and $Y$ are locally compact separable metric spaces, $\mu\in \M_+^1(X)$,  and $f$ : $X\rightarrow Y$ is a Borel map we shall denote by $f\sharp \mu$ the push forward of $\mu$ through $f$ i.e. the element of $\M_+^1(Y)$ defined by $f\sharp \mu (B)=\mu (f^{-1}(B))$ for every Borel subset $B$ of $Y$. 

\smallskip
In the sequel, $\Leb^d$ denotes the $d$-dimensional Lebesgue measure. If $\mu$ and $\nu$ are in $\M_+^1(\R^d)$ then $\frac{d\mu}{d\nu}$ denotes the Radon-Nikodym derivative of $\mu$ with respect to $\nu$. We shall write $\mu<<\nu$ to express that $\mu$ is absolutely continuous with respect to $\nu$, in which case, slightly abusing notations, we will identify $\mu$ with the Radon-Nikodym derivative $\frac{d\mu}{d\nu}$.

\smallskip

The data of our problem are $\Omega$ (its closure $\omb$ modelling the city, say) which is some open bounded  convex subset of $\R^2$, two probability measures, $\mu_0$ and $\mu_1$ in $\M_+^1(\omb)$, giving respectively the distribution of residents and services in the city $\omb$.  The set of transportation plans associated to $\mu_0$ and $\mu_1$ will be denoted $\Pi(\mu_0,\mu_1)$ it consists of the probability measures on $\omb\times \omb$ having $\mu_0$ and $\mu_1$ as marginals:
\begin{equation}\label{deftp}
\Pi(\mu_0,\mu_1):=\{\gamma \in \M_+^1(\omb\times \omb) \mbox{ : } \pi_0\sharp \gamma=\mu_0,\; \pi_1\sharp \gamma=\mu_1\}
\end{equation} 
where $(\pi_0(x,y),\pi_1(x,y)):=(x,y)$, stand for the canonical projections ($x$ and $y$ in $\omb$).

Introducing congestion naturally leads to consider spaces of paths, lengths of such paths and sets of probability measures on sets of paths. From now, on we shall denote:

\begin{itemize}
\item $C:=W^{1,\infty}([0,1], \omb)$, viewed as a subset of $C^0([0,1],\R^2)$,
\item $C^{x, y}:=\{\sigma\in C:\ \sigma(0)=x,\ \sigma(1)=y\}$ ($x,\ y$ in  $\omb$),
\item $l(\sigma):=\int_0^1 \vert \dot\sigma(t)\vert\ dt$, the length of $\sigma\in C$,
\item for $\sigma\in C$, $\tilsig$ denotes the arclength reparameterization of $\sigma$ belonging to $C$, hence $\vert \dot{\tilsig}(t)\vert=l(\sigma)=l(\tilsig)$ for a.e. $t\in [0,1]$,  
\item $\tilC:=\{\sigma \in C \; :\; \vert \dot{\sigma}\vert \mbox{ is constant}\}=\{\tilsig, \; \sigma\in C\}$, 
\item for $Q\in \M_+^1(C)$, we define $\tilQ\in \M_+^1(\tilC)$ as the push forward of $Q$ through the map $\sigma\mapsto \tilsig$,
\item for $\varphi\in C^0(\omb,\R)$ and $\sigma\in C$, we define 
\[L_{\varphi}(\sigma):=\int_0^1 \varphi(\sigma(t))\vert \dot{\sigma}(t)\vert dt=l(\sigma)\int_0^1 \varphi(\tilsig(t))dt,\] 
\item $e_0(\sigma):=\sigma(0)$, $e_1(\sigma):=\sigma(1)$,  for all $\sigma\in C^0([0,1],\R^2)$.
\end{itemize}

\subsection{Traffic congestion modelling}\label{tcm}

The classical Monge-Kantorovich optimal transportation problem  for a given cost function $c\in C^0(\omb\times\omb,\R)$ is:
\begin{equation}\label{mk}
\inf \left\{  \int_{\omb\times \omb} c(x,y) d\gamma(x,y) \; :\; \gamma\in \Pi(\mu_0,\mu_1) \right\}. 
\end{equation}
Note that, in the linear problem \pref{mk}, the cost of transporting one unit of mass from $x$ to $y$, $c(x,y)$, is given and does not depend on the path(s) followed by the mass from $x$ to $y$. In order to take into account congestion effects, we explicitely introduce probabilities over $C^{x,y}$ as part of the optimization problem. More precisely, the overall transportation cost will depend not only on the transportation plan $\gamma\in\Pi(\mu_0,\mu_1)$ but also on the way travelers  commuting from $x$ to $y$ use the different possible paths $\sigma\in C^{x,y}$. In the sequel,  the way commuters from $x$ to $y$ are split according to the different paths will be given  by a  probability measure $p^{x,y}$ on $C^{x,y}$. Put differently, $p^{x,y}(\Sigma)$ is the proportion of travelers from $x$ to $y$ using a path $\sigma\in \Sigma\subset C^{x,y}$. This naturally leads to the following definition:

\begin{defi}
A transportation strategy consists of a pair $(\gamma, p)$ with $\gamma\in \Pi(\mu_0,\mu_1)$ 	and where $p=(p^{x,y})_{(x,y)\in\omb\times\omb}$ is a Borel family of probability measures on $C$ such that $p^{x,y}(C^{x,y})=1$ for $\gamma$-a.e. $(x,y)\in\omb\times\omb$.  
\end{defi}  

There results, from the use of a transportation strategy $(\gamma,p)$, an overall \emph{traffic intensity} $I_{\gamma,p}\in \M_+(\omb)$ defined by
\begin{eqnarray}\label{defiIgammap}
\int_{\omb}  \varphi(x) d I_{\gamma,p}(x):=\int_{\omb\times \omb}  \left( \int_{C^{x,y}} \left( \int_0^1 \varphi(\sigma(t))\vert \dot{\sigma}(t)\vert  dt \right) dp^{x,y}(\sigma)  \right) d\gamma(x,y)\nonumber\\
=\int_{\omb\times \omb}  \left( \int_{C^{x,y}} L_{\varphi}(\sigma) dp^{x,y}(\sigma)  \right) d\gamma(x,y), \quad \forall \varphi \in C^{0}(\omb,\R)\nonumber\\
 \label{deftrafintgp}
\end{eqnarray} 
and an overall probability over paths $Q_{\gamma,p}\in \M_+^1(C)$ given by $Q_{\gamma,p}=p^{x,y}\otimes \gamma$, i.e.:
\begin{equation}\label{defQ}
\int_C F(\sigma) dQ_{\gamma,p}(\sigma)=\int_{\omb\times \omb}  \left( \int_{C^{x,y}} F(\sigma) dp^{x,y}(\sigma)  \right) d\gamma(x,y)\quad \forall F \in C^{0}(C,\R).
\end{equation}
One could consider the probability $Q_{\gamma,p}$ as if it represented the total number of travelers that use a path $\sigma\in \Sigma$ given the global transportation strategy $(\gamma,p)$. 

\smallskip

 Let us remark that if we set $Q:=Q_{\gamma,p}\in \M_+^1(C)$ then $I_{\gamma,p}$ only depends on $Q$, and can be written as $I_{\gamma,p}=i_Q\in  \M_+(\omb)$ where $i_Q$ is defined for every $Q\in \M_+^1(C)$ by:
\begin{equation}\label{defiQ}
\int_{\omb} \varphi(x) d i_Q(x)=\int_C L_{\varphi}(\sigma) d Q(\sigma),\quad \forall \varphi \in C^{0}(\omb,\R).
\end{equation}
Let us also remark that since $L_{\varphi}(\sigma)=L_{\varphi}(\tilsig)$ one has $i_Q=i_{\tilQ}$, for all $Q\in \M_+^1(C)$. Finally, let us note that the total mass of $i_Q$ is the average length with respect to $Q$:
\begin{equation}\label{toatalmassiq}
i_Q(\omb)=\int_C l(\sigma)dQ(\sigma).
\end{equation} 
 if the probability $Q$ is concentrated on injective curves, one could also express the measure $i_Q$ through $\huno-$integrals. In this same case, namely 
\begin{rem}\label{huno}
If a curve $\sigma$ is injective, then one could also write $l(\sigma)=\huno(\sigma([0,1]))$ and $L_\varphi(\sigma)=\int_{\sigma([0,1])}\varphi\,d\huno$. Moreover, if for $\gamma-$a.e. $(x,y)$ the probability $p^{x,y}$ is concentrated on the set of injectives curves from $x$ to $y$, one could also define the measure $I_{\gamma,p}$ by replacing the integral with respect to $\vert \dot{\sigma}(t)\vert  dt $ in \eqref{defiIgammap} with an integral in $d\huno$. Notice moreover that for every  Borel subset $A\subset\omb$ one would have:
\[\begin{split}
I_{\gamma,p}(A)=\int_{\omb\times \omb}  \left( \int_{C^{x,y}}\! \huno(A\cap \sigma) dp^{x,y}(\sigma)\right) d\gamma(x,y)= \int_C \!\huno(A\cap \sigma) dQ_{\gamma,p}(\sigma).
\end{split}\]
If we imagine that for each $\sigma\in C^{x,y}$, the mass of travelers commuting on $\sigma$ is uniformly distributed on $\sigma$, this means that $I_{\gamma,p}(A)$ represents the cumulative traffic through the region $A$. The same formula stays true, under no injectivity assumption, if we replace $\huno(A\cap\sigma)$ with $L_{I_A}(\sigma)$ and in this case the cumulative traffic takes into account the number of times a path $\sigma$ passes through the points of $A$.
\end{rem}

In the sequel, it will be convenient to formulate our optimization problem in terms of $Q=Q_{\gamma,p}$ rather than in the transportation strategy $(\gamma,p)$. To that end, we shall use the following:

\begin{lem}\label{qq}
Let us define
\[\QQ:=\{Q_{\gamma, p} \; :\; (\gamma,p) \mbox{ transportation strategy}\}\]
then one has
\[\QQ=\{Q\in\M_+^1(C) \; : \; e_0\sharp Q= \mu_0, \; e_1\sharp Q= \mu_1\}.\]
\end{lem}

\begin{proof}
If $(\gamma,p)$ is a transportation strategy then $e_0\sharp Q_{\gamma,p}= \pi_0\sharp \gamma=\mu_0$, and $e_1\sharp Q_{\gamma,p}=\pi_1\sharp \gamma= \mu_1$. Now let $Q\in\M_+^1(C)$ be such that $e_0\sharp Q= \mu_0$, $e_1\sharp Q= \mu_1$. If we define $\gamma:=(e_0,e_1)\sharp Q$, we have $\gamma\in \Pi(\mu_0,\mu_1)$. It then follows from the disintegration theorem (see \cite{dm}) that  there exists $p=(p^{x,y})_{(x,y)\in\omb\times\omb}$  a Borel family of probability measures on $C$ such that $p^{x,y}(C^{x,y})=1$ for $\gamma$-a.e. $(x,y)\in\omb\times\omb$ and $Q= p^{x,y}\otimes \gamma$. Hence  $Q=Q_{\gamma,p}$ for a transportation strategy $(\gamma,p)$.
\end{proof}

At this point, a natural way to model traffic congestion is, for a given transportation strategy $(\gamma,p)$, to consider that the transportation cost per unit of mass between  $x$ and $y$  is given by
\begin{equation}\label{defcgp}
c_{\gamma, p} (x,y)=\int_{C^{x,y}} L_{G_{I_{\gamma,p}}}(\sigma) dp^{x,y}(\sigma)
\end{equation}
where $G_{I_{\gamma,p}}$ is a nonnegative function which depends (in a way that will be  specified later on) on the traffic intensity $I_{\gamma,p}$.  The optimal transportation with traffic congestion then takes the form (to be compared with the usual Monge-Kantorovich problem \pref{mk}):
\begin{equation}\label{pbmefi}
\inf  \left\{  \int_{\omb\times \omb}  c_{\gamma,p} (x,y) d\gamma(x,y) \;  : \;  (\gamma,p) \mbox{ transportation strategy} \right\}.
\end{equation}
Setting $Q=Q_{\gamma,p}$ and using formally \pref{defiQ}, we see that the total transportation cost in \pref{pbmefi} can be rewritten as:
\[\int_{\omb\times \omb}  c_{\gamma,p} (x,y) d\gamma(x,y) =\int_{C} L_{G_{i_Q}}(\sigma) dQ(\sigma)
=\int_{\omb} G_{i_Q}(x) d i_Q(x).\]
Hence using lemma \ref{qq}, we can reformulate \pref{pbmefi} in terms of $Q$ only:
\begin{equation}\label{pbmenQ}
\inf \left \{  \int_{\omb} G_{i_Q}(x) d i_Q(x )  \; :\; Q\in\QQ \right\}
\end{equation}
Note that in the definition \pref{defcgp}, it is required that $G_{I_{\gamma,p}}$ is  continuous (or at least l.s.c.) whereas the form \pref{pbmenQ} allows for more general forms of congestion through $i\mapsto G_i$. From now on, we assume that  $G$ has the following local form:
\begin{equation}
G_i(x)=g\left(\frac{di}{d \Leb^2}(x)\right) ,
\end{equation}
where $\frac{di}{d \Leb^2}$ is the radon-Nicodym derivative of $i$ with respect to the Lebesgue measure and $g$ is a nondecreasing function $\R_+\rightarrow \R_+$ such that the function $H$  defined by $H(z)=zg(z)$ for all $z\in \R_+$ is convex and superlinear (i.e. $\lim_{z\to +\infty}g(z)=+\infty$).

\smallskip

The optimization problem we shall study now then reads as: 
\begin{equation}\label{lepbme0}
\inf_{Q\in \QQ} \HH(i_Q) \mbox{ where }
\ \HH(i)=\left\{
\begin{array}{l}
\int_{\Omega}H(i(x))dx  \mbox{ if $i<< \Leb^2$,}
\\ +\infty
\mbox{ otherwise.}\end{array}\right.
\end{equation}
In the sequel, we shall say that a transportation strategy $(\gamma, p)$ is optimal if $Q_{\gamma,p}$ solves \pref{lepbme}.

\begin{rem}
It will be clear in the sequel that the probability $Q_{\gamma,p}$ associated to an optimal transportation strategy $(\gamma,p)$ will be concentrated on injective curves, so that the interpretation in terms of $\huno-$integrals (see Remark \ref{huno}) may apply.
\end{rem}

\subsection{Existence of minimizers}

From now on, we make the following assumptions:
\begin{itemize}
\item $H$ is convex and nondecreasing on $\R_+$ with $H(0)=0$,
\item there exists $q>1$, and positive constants $a$ and $b$ such that\\
 $az^q\leq H(z) \leq b( z^q+1)$ for all $z\in \R_+$, 
\item $H$ is differentiable on $\R_+$, and there exists a positive constant $c$ such that $0\leq H'(z)\leq c(z^{q-1}+1)$,  for all $z\in \R_+$, 
\item the following set
\begin{equation}
\QQq:=\{Q\in \QQ \; : \; i_Q\in L^q\}
\end{equation}
 is nonempty.
\end{itemize}

These assumptions enable us to simply rewrite \pref{lepbme0} as:
\begin{equation}\label{lepbme}
\inf_{Q\in\QQq} \int_{\Omega} H(i_Q(x))dx.\end{equation}

\begin{rem}\label{domainenonvide}
Let us discuss the assumption that $\QQq\neq \emptyset$ which, at first glance,  may seem difficult to check. In order to have the existence of a $Q\in \QQ$ such that $i_Q\in L^q$ it is sufficient that  $\mu_0$ and $\mu_1$ are in $L^q$. This result, which is not obvious, follows from the regularity  results of De Pascale and Pratelli (see \cite{dpp} and \cite{dpp2}) who proved that $L^q$ regularity of $\mu_0$ and $\mu_1$ implies that for $\gamma$ solving the Monge-Kantorovich problem \pref{mk} with $c(x,y)=\vert x-y\vert$ and $p^{x,y}=\delta_{[x,y]}$ (the Dirac mass at the segment $[x,y]$) for every $x$ and $y$ the corresponding traffic density $I_{\gamma,p}$ is $L^q$.  It is not necessary however that $\mu_0$ and $\mu_1$ are absolutely continuous for the assumption to be satisfied: let us consider for instance the case where  $\omb=[0,1]^2$ and $\mu_0$ and $\mu_1$ are respectively the one-dimensional Hausdorf measures of  the segments $[(0,0), (0,1)]$ and $[(1,0),(1,1)]$. If we define $\gamma:=(\mathrm{id}, \mathrm{id}+(1,0))\sharp\mu_0$   and  $p^{x,y}=\delta_{[x,y]}$  then a straightforward computation shows that $I_{\gamma,p}$ is uniform on $[0,1]^2$.
\end{rem}

\smallskip

Under the assumptions above, we are going prove that \pref{lepbme} admits a solution.  The proof of existence involves some preliminary lemmas.

\begin{lem}\label{lsc}
For any   $\varphi \in C^0(\omb,\R_+)$, 
$L_{\varphi}$ is lower semi-continuous on $C$ for the uniform topology, indeed for any $\sigma\in C$, one has:
\begin{multline}\label{LPHI}
L_{\varphi}(\sigma)=\sup\Big\{\sum_{i=1}^n \left(\inf_{[t_i,t_{i+1}]}(\varphi\circ\sigma)\right)\vert \sigma(t_{i+1})-\sigma(t_i)\vert:\\ 
([t_i, t_{i+1}])_i \mbox{ is a subdivision of }[0,1]\Big\}. 
\end{multline} 

\end{lem}

\begin{proof}
For any subdivision $([t_i, t_{i+1}])_{i=1,...n}$, we have:
\begin{eqnarray*}  
L_{\varphi}(\sigma)&=&\sum_{i=1}^n \int_{t_i}^{t_{i+1}} \varphi(\sigma(t))\vert\dot \sigma(t)\vert\ dt\\
&\ge& \sum_{i=1}^n \inf_{[t_i,t_{i+1}]} (\varphi\circ\sigma) \int_{t_i}^{t_{i+1}} \vert \dot\sigma(t)\vert\ dt\\
&\ge& \sum_{i=1}^n \inf_{[t_i,t_{i+1}]}(\varphi\circ\sigma)\vert \sigma(t_{i+1})-\sigma(t_i)\vert.
\end{eqnarray*} 
Taking the supremum over all such divisions, we get:
\begin{equation*}
L_{\varphi}(\sigma)\ge\sup\Big\{\sum_{i=1}^n \inf_{[t_i,t_{i+1}]}(\varphi\circ\sigma)\vert \sigma(t_{i+1})-\sigma(t_i)\vert:\\ 
([t_i, t_{i+1}])_{i} \mbox{ is a subdivision of }[0,1]\Big\}. 
\end{equation*} 
Let us prove the converse inequality. Let $\e>0$, since $\varphi\circ\sigma$ is uniformly continuous, there is a $\delta>0$ such that:
\begin{equation*}
\forall t,t'\in[0,1]^2,\quad \left(\vert t-t'\vert\le\delta\Rightarrow \vert \varphi(\sigma(t)) -\varphi(\sigma(t'))\vert\le \e\right).
\end{equation*}
For any subdivision $([t_i, t_{i+1}])_{i=1,...n}$ such that $\vert t_i-t_{i+1}\vert\le \delta$ for all $i$, we have:
\begin{eqnarray*}  
&&L_{\varphi}(\sigma)
\le \sum_{i=1}^n (\inf_{[t_i,t_{i+1}]} (\varphi\circ\sigma) +\e) \int_{t_i}^{t_{i+1}} \vert \dot\sigma(t)\vert\ dt\\
&&= \sum_{i=1}^n (\inf_{[t_i,t_{i+1}]} (\varphi\circ\sigma) +\e) 
\sup\big\{ \sum_j \vert \sigma(\tau_j)-\sigma(\tau_{j+1})\vert:\\
&&\hspace{3.5cm} ([\tau_j, \tau_{j+1}])_j\mbox{ is a subdivision of }[t_i,t_{i+1}]\big\}\\
&&\le\sup\big\{  \sum_i\sum_j(\inf_{[\tau_j,\tau_{j+1}]} (\varphi\circ\sigma) +\e) \vert \sigma(\tau_j)-\sigma(\tau_{j+1})\vert:\\
&&\hspace{3.5cm}([\tau_j, \tau_{j+1}])_j 
\mbox{ is a subdivision of }[t_i,t_{i+1}]\big\}\\
&&= \sup\Big\{\sum_{i=1}^n (\inf_{t\in [t_i,t_{i+1}]}(\varphi\circ\sigma)+\e)\vert \sigma(t_{i+1})-\sigma(t_i)\vert:\\
&&\hspace{3.5cm}
([t_i, t_{i+1}])_i \mbox{ is a subdivision of }[0,1]\Big\}.
\end{eqnarray*} 
As this last inequality is true for any $\e>0$ we get (\ref{LPHI}). The lower semi-continuity is then obvious since, by (\ref{LPHI}), 
$L_{\varphi}$ is the supremum of family of  lower semi-continuous functions on $C^0([0,1], \omb)$.
\end{proof}

\begin{lem}\label{tight}
Let $(Q_n)_n\in \M_+^1(C^0([0,1],\R^2))^{\N}$ such that $Q_n(C)=1$ for all $n$ and there exists a constant $M>0$ such that:
$$\sup_n \int_C l(\sigma)\ dQ_{n}(\sigma)\le M.$$ 
Then the sequence $(\tilQ_n)_n$ is tight and admits a subsequence that converges weakly $*$ to a probability $Q$ such that $Q(C)=1$. 
\end{lem}
\begin{proof}
The tightness of $(\tilQ_n)_n$ easily follows from the inequality:
\begin{eqnarray}
\tilQ_n \left(\{\sigma\in \tilC:\ \vert \dot\sigma\vert> K\}\right)&= & 
Q_n\left(\left\{\sigma\in C:\ l(\sigma)> K\right\}\right)\nonumber\\
&\le& {1\over K}   \int_C l(\sigma)\ dQ_{n}(\sigma)\label{estimtight}.
\end{eqnarray}
By Prokhorov theorem, we may therefore assume, passing to a subsequence if necesseary, that  $(\tilQ_n)_n$ converges weakly $*$ to $Q\in \M_+^1(C^0([0,1],\R^2)$.   It remains to show that  $Q(C)=1$. For $K>0$ let us define  $C_K:=\{\sigma \in C \; :\;  \vert \dot{\sigma}\vert \leq K\}$, then Inequality \pref{estimtight} and the fact that the measures $\tilQ_n$ are concentrated on $\tilC$ yield
\[\sup_n \,\tilQ_n(C\backslash C_K)= \sup_n \,\tilQ_n(\tilC\backslash C_K)\le {M\over K},\]
for every $K>0$, which implies
\begin{eqnarray*}
1=\limsup_n \tilQ_n(C)&\leq & \limsup_n \tilQ_n(C_K) +\limsup_n \tilQ_n(C\setminus C_K) \\
&\le& Q(C_K)+ {M\over K}.
\end{eqnarray*}
Letting $K$ tend to $\infty$, we then get $Q(C)=\sup_K Q(C_K)=1$. 
\end{proof}

\begin{lem}\label{I-IQ}
Let $(Q_n)_n$ be a sequence in $\M_+^1(C)$ that converges weakly $*$ to some $Q\in  \M_+^1(C)$. If there exists $i\in \M_+(\omb)$ such that $i_{Q_n}$ converges weakly $*$ to $i$ in  $\M_+(\omb)$
 then we have $i_Q\le i$.  
\end{lem}

\begin{proof}
Let $\varphi\in C^0(\omb, \R_+)$, we first have:
\[\int_{\omb}\varphi d i=\lim_{n} \int_{\omb}\varphi d  i_{Q_n} =\lim_{n} \int_C   L_{\varphi} dQ_n\]
it easily follows from lemma \ref{lsc} that $Q\mapsto \int_C   L_{\varphi} dQ$ is l.s.c. for the weak $*$ topology of $\M_+^1(C)$, we then have:
\[\int_{\omb}\varphi d i
\ge\int_C L_{\varphi}dQ=\int_{\omb} \varphi di_Q.\qedhere\] 
\end{proof}

Now, we are in position to prove:

\begin{thm}\label{existence}
The minimization problem \pref{lepbme} admits a solution.
\end{thm}

\begin{proof}
Our assumptions imply that the value of \pref{lepbme} is finite.  Let $(Q_n)_n$ be some minimizing sequence of  \pref{lepbme}. From the identity $i_{Q}=i_{\tilQ}$, we may assume $Q_n=\tilQ_n$ for all $n$.  We deduce from our growth condition on $H$, that $(i_{Q_n})_n$ is bounded in $L^q$. On the one hand, extracting a subsequence if necessary, we may therefore assume that $(i_{Q_n})_n$ converges weakly in $L^q$ to some $i$. On the other hand, since $i_{Q_n}$ is bounded in $L^q$ and hence in $L^1$ we have
\[\sup_n \int_C l(\sigma)dQ_n(\sigma)=\sup_n\int_{\Omega} i_{Q_n}<+\infty.\]
Moreover $Q_n=\tilQ_n$ and we deduce from lemma \ref{tight} that (up to some subsequence)  $(Q_n)_n$ weakly $*$ converges to some $Q$ in $\M_+^1(C)$. Since $\QQ$ is  obviously weakly $*$ closed, we have $Q\in \QQ$ and lemma \ref{I-IQ} implies that $i_Q\leq i$ (consequently to this enaquality $i_Q$ is absolutely continuous). From the monotonicity and convexity of $H$  we then have: 
\[\int_{\Omega} H(i_Q(x))dx \leq \int_{\Omega} H(i(x))dx \leq \liminf_n \int_{\Omega} H(i_{Q_n}(x))dx\]
which proves that $Q$ solves \pref{lepbme}.  
\end{proof}

Let us remark that if $H$ is furthermore assumed to be strictly convex then if $Q_1$ and $Q_2$ solves \pref{lepbme} then $i_{Q_1}=i_{Q_2}$ so that the optimal traffic intensity is unique (of course, this does not imply in general that $Q_1=Q_2$ or that the corresponding optimal transportation strategy is unique).

\section{Characterization of the minimizers}

In the sequel, we shall denote by $q^*$ the conjugate exponent of $q$, given by $q^*=q/(q-1)$. 

\subsection{Optimality conditions}\label{foc}

The variational inequalities characterizing solutions of the convex problem \pref{lepbme} can be expressed as follows

\begin{prop}
$\oQ\in \QQq$  solves \pref{lepbme} if and only if
\begin{equation}\label{varineq}
\int_{\Omega} \oxi i_{\oQ}=\inf \left\{ \int_{\Omega} \oxi i_{Q} \; :\; Q\in\QQq\right\} \mbox{ with } \oxi:=H'(i_ {\oQ})\in L^{q^*}.
 \end{equation}
\end{prop}
\begin{proof}
Assume that  $\oQ$ solves \pref{lepbme}, then for every $Q\in\QQq$, one has:
\[\begin{split}
0  \leq  &\lim_{\e\rightarrow 0^+} \frac{1}{\eps} [\HH(i_{\oQ +\e(Q-\oQ)})-\HH(i_{\oQ})]
= \lim_{\e\rightarrow 0^+} \frac{1}{\eps} [\HH(i_{\oQ} +\e (i_{Q}-i_{\oQ}))-\HH(i_{\oQ})]\\
=&\int_{\Omega} H'(i_{\oQ})(i_{Q}-i_{\oQ})= \int_{\Omega} \oxi (i_{Q}-i_{\oQ})
\end{split}\]
which proves \pref{varineq}. Conversely, if $\oQ\in\QQq$ satisfies \pref{varineq}, then by convexity of $H$, for every  every $Q\in\QQq$, one has:
\[\HH(i_Q)-\HH(i_{\oQ})\geq \int_{\Omega} \oxi (i_Q-i_{\oQ})\geq 0.\qedhere\]
\end{proof}

The next paragraphs will be devoted to investigate the precise meaning of \pref{varineq}.  Before going further, let us do some formal manipulations to give a formal interpretation of  \pref{varineq}  in terms of optimal transportation strategy. Let us assume that $\oQ$ solves \pref{lepbme} and let us write $\oQ=Q_{\ogm,\ovp}$ for some (optimal) transportation strategy $(\ogm,\ovp)$ and define $\oxi:=H'(i_{\oQ})$, then \pref{varineq} formally can be rewritten as:
\[\begin{split}
\int_{\Omega} \oxi i_{\oQ} = &\int_C L_{\oxi} (\sigma)d\oQ(\sigma)\\
=&\int_{\omb\times \omb} \left( \int_{C^{x,y}} L_{\oxi}(\sigma)d\op(\sigma)\right)d\ogm(x,y)\\
=&\inf_{(\gamma, p) \mbox{ transp. strategy} } \int_{\omb\times \omb} \left(\int_{C^{x,y}} L_{\oxi}(\sigma)dp^{x,y}(\sigma)\right)d\gamma(x,y)\\
=&\inf_{\gamma\in \Pi(\mu_0,\mu_1)} \int_{\omb\times \omb} \left(\inf_{p\in\M_+^1(C^{x,y})} \int_{C^{x,y}} L_{\oxi}(\sigma)dp(\sigma)\right)d\gamma(x,y)\\
=&\inf_{\gamma\in \Pi(\mu_0,\mu_1)} \int_{\omb\times \omb} \left(\inf_{\sigma\in C^{x,y}}  L_{\oxi}(\sigma)\right)d\gamma(x,y)
\end{split}\]
defining (again formally) the transportation cost:
\[c_{\oxi}(x,y)=\inf_{\sigma\in C^{x,y}}  L_{\oxi}(\sigma),\]
we then firstly have:
\[\int_{\omb\times \omb} c_{\oxi}(x,y)d\ogm(x,y)\leq \int_{C} L_{\oxi} d\oQ=  \inf_{\gamma\in \Pi(\mu_0,\mu_1)} \int_{\omb\times \omb} c_{\oxi}(x,y)d\gamma(x,y) \] 
so that $\ogm$ solves the Monge-Kantorovich problem:
\[\inf_{\gamma\in \Pi(\mu_0,\mu_1)} \int_{\omb\times \omb} c_{\oxi}(x,y)d\gamma(x,y).\]
Secondly:
\[\begin{split}
\int_{C} L_{\oxi}(\sigma) d\oQ(\sigma)= \int_{\omb\times \omb} c_{\oxi}(x,y)d\ogm(x,y)\\
=\int_C  c_{\oxi}(\sigma(0),\sigma(1))d\oQ(\sigma)\end{split}\]
and since $L_{\oxi}(\sigma)\geq c_{\oxi}(\sigma(0),\sigma(1))$, we get
\[L_{\oxi}(\sigma)=c_{\oxi}(\sigma(0),\sigma(1)) \quad\mbox{ for }  \oQ \mbox{-a.e. } \sigma.\]
or, in an equivalent way, for $\ogm$-a.e. $(x,y)$ one has:
\[ L_{\oxi}(\sigma)=c_{\xi}(x,y)\; \quad\mbox{ for }  \op \mbox{-a.e. } \sigma.\]
Since $\oxi$ is only $L^{q^*}$, $L_{\xi}$ and $c_{\oxi}$ are not well-defined and the previous arguments are purely formal.  In paragraph \ref{cxinr}, we will extend the definition of $c_{\xi}$ to the case where $\xi$ is  only $L^{q^*}$ under the additional assumption $q<2$. This will enable us to make the formal argument above rigorous and to characterize optimal transportation strategies in paragraph \ref{ots}.

\subsection{The transportation cost $\oc_{\xi}$  when $\xi$ is  $L^{q^*}$}\label{cxinr}

For a non-negative function $\xi \in C^0(\Omega)$ we define
$$c_\xi(x,y)=\inf\{L_{\xi}(\sigma)\,:\,\sigma\in C^{x,y}\}.$$

\begin{prop}\label{cxicomp}
Let us assume that $q<2$ and define $\alpha:=1-2/q^*$, then there exists a non-negative constant $C$ such that for every $\xi\in C^0(\omb,\R^+)$ and every $(x_1,y_1,x_2,y_2)\in\Omega^4$, one has:
\begin{equation}\label{holderest}
\vert c_{\xi}(x_1,y_1)-c_{\xi}(x_2,y_2)\vert \leq C\Vert \xi\Vert_{L^{q^*}(\Omega)} \left( \vert x_1-x_2\vert^{\alpha}+  \vert y_1-y_2\vert^{\alpha} \right).
\end{equation}
Consequently, if $(\xi_n)_n\in  C^0(\omb,\R^+)^{\N}$ is bounded in $L^{q^*}$, then $(c_{\xi_n})_n$ admits a subsequence that converges in $C^0(\omb\times \omb,\R_+)$. 
\end{prop}

\begin{proof}
 Let $\xi\in C^0(\omb,\R^+)$ and $x,\ y\in \Omega^2$. For  $k>0$ let $\sigma_k\in C^{x,y}$ be such that
\[\int_0^1 \xi(\sigma_k(t))\vert \dot \sigma_k(t)\vert\ dt\leq c_{\xi}(x,y)+ {1\over k}.\]
Then for all $\e>0$ and $z$ such that $y+\e z\in \Omega$ and $t_0\in (0,1)$ we consider the following element 
of $C^{x, y+\e z}$: 
\begin{equation*}
\sigma_{k,t_0}(t):=\left\{
\begin{array}{c l}
\sigma_k\left({t\over t_0}\right) &\mbox{ if }t\in[0,t_0]\\  
y+ \left({t-t_0\over 1-t_0}\right) \e z &\mbox{ if }t\in[t_0,1].\\
\end{array}
\right.
\end{equation*}
We then  have, for all $k>0$:
\begin{eqnarray*}
c_{\xi}(x,y+\e z)&\le& \int_0^1 \xi(\sigma_{k,t_0}(t))\vert \dot \sigma_{k,t_0}(t)\vert \ dt\\ 
&=&\int_0^{1} \xi (\sigma_k(t))\vert \dot \sigma_k(t)\vert \ dt 
+ \int_{0}^1 \xi (y+ t \e z )\e \vert z\vert \ dt\\
&\le & c_{\xi}(x,y)  +  \e \vert z\vert\int_{0}^1 \xi (y+ t \e z ) dt +{1\over k}.
\end{eqnarray*}
Now we let $k$ tend to $\infty$ and we get
\[{c_{\xi}(x,y+\e z)-c_{\xi}(x,y)\over \e}\le \vert z\vert\int_0^1 \xi(y+t\e z)\ dt,\]
and by a similar argument
\[{c_{\xi}(x,y)-c_{\xi}(x,y+\e z)\over \e}\le \vert z\vert\int_0^1 \xi(y+(1-t)\e z)\ dt.\]
This implies that $c_{\xi}(x,.)\in W^{1,\infty}$ and:
\begin{equation}\label{EST}
\vert \nabla_y c_{\xi}(x,.)\vert\le \vert \xi(.)\vert, \; \mbox{ for all} \;  x.
\end{equation}
By symmetry we also have
\begin{equation}\label{ESTx}
\vert \nabla_x c_{\xi}(.,y)\vert\le \vert \xi(.)\vert, \; \mbox{ for all} y.
\end{equation}
Since $q^*>2$, we deduce from \pref{EST}, \pref{ESTx} and Morrey's Theorem (see \cite{brezis}, Chapter IX), that  there is a constant $C$ such that:
\[\begin{split}
\vert c_{\xi}(x,y_1)-c_{\xi}(x,y_2)\vert \leq C\Vert \xi\Vert_{L^{q^*}}  \vert y_1-y_2\vert^{\alpha} 
,\;   \mbox{for all }  x, y_1, y_2 \mbox { in } \Omega,\\
\vert c_{\xi}(x_1,y)-c_{\xi}(x_2,y)\vert \leq C\Vert \xi\Vert_{L^{q^*}}  \vert x_1-x_2\vert^{\alpha} 
,\;   \mbox{for all }  x_1, x_2, y \mbox { in } \Omega.
\end{split}\]
This proves \pref{holderest}. The second claim in the proposition then follows from \pref{holderest}, the identity $c_{\xi_n}(x,x)=0$ and Ascoli's theorem. 
\end{proof}

From now on, we further assume that $q<2$. For a non-negative function $\xi\in L^{q^*}(\Omega)$ we then define
$$\overline{c}_\xi(x,y)=\sup\left\{c(x,y)\,:c\in\mathcal{A}(\xi)\right\},$$
where 
$$\mathcal{A}(\xi)=\left\{\lim_n c_{\xi_n}\,\mbox{ in }C^0(\omb\times\omb)\,:\,(\xi_n)_n\in C^0(\omb),\,\xi_n\geq 0,\,
\xi_n\to\xi\,\mbox{ in }L^{q^*}\right\}.$$

\begin{rem}\label{cxifaible}
The definition of $\oc_{\xi}$ is unchanged if one replaces $\xi_n\to\xi\,\mbox{ in }L^{q^*}$ by $\xi_n\weakto \xi\,\mbox{ in }L^{q^*}$ in the definition of $\A(\xi)$. Indeed, if we do so, we obviously obtain a function which is larger than $\oc_{\xi}$. Now, let us assume that $\xi_n\weakto \xi\,\mbox{ in }L^{q^*}$, and $c_{\xi_n}$ converges to $c$ in $C^0(\omb\times\omb)$, using Mazur's Lemma there exists a sequence $\eta_n$ which converges strongly to $\xi$ and such that each $\eta_n$ is in the convex hull of $\{\xi_k, \; k\leq n\}$. It is clear that for fixed $x$, $y$, $\xi\to c_{\xi}(x,y)$ is concave hence $c(x,y)=\lim c_{\xi_n}(x,y) \leq \limsup  \; c_{\eta_n}(x,y) \leq \oc_{\xi}(x,y)$.
\end{rem}

When $\xi$ is continuous, one has:

\begin{lem}\label{ccoincid}
If $\xi$ is continuous and non-negative, then $\overline{c}_\xi=c_\xi$.
\end{lem}
\begin{proof}

 The inequality $\overline{c}_\xi\geq c_\xi$ is obvious, 
as one can always choose the constant sequence $\xi_n=\xi$ in the 
definition of $\overline{c_\xi}$. Let us show now the opposite inequality. Take $x,y\in\Omega$, $\ve>0$ 
and $\sigma\in C^{x,y}$ such that $L_{\xi}(\sigma)<c_\xi(x,y)+1/k$. We can choose $\sigma$ so that it is 
piecewise linear, by density of this kind of curves and using the continuity of $\xi$. Let $(S_i)_{i=0,\dots,m-1}$ 
be the segments which compose $\sigma$ with $S_i=x_i x_{i+1}$, $x_0=x$ and $x_m=y$. Let us fix, moreover, 
a sequence $\xi_n\to\xi$ such that $c_{\xi_n}\to c$.  Now, we want to prove $c\leq c_{\xi}$. Fix a small number $\delta>0$ and for any $\alpha\in [0,\delta]$ let us define a curve $\sigma^\alpha$ in the following way: 
let $R$ be the clockwise $90$ degrees rotation in the plane;
let $x'_i(\alpha)$ and $x''_i(\alpha)$ be the only points such that $x_ix'_i(\alpha)=\alpha Re_i$ and $x_{i+1}x''_i(\alpha)=\alpha Re_i$,
where $e_i$ is the tangent unit vector to $\sigma$ in the $S_i$ part; define $\sigma^\alpha$ by linking any point $x'_i(\alpha)$ to 
$x''_i(\alpha)$ by some segments $S'_i(\alpha)$ and $x''_i(\alpha)$ to $x'_{i+1}(\alpha)$ by some arcs $A_{i+1}(\alpha)$ with center 
$x_{i+1}$ and radius $\alpha$. 
In this way we have $\sigma^\alpha\in C^{x_\alpha,y_\alpha}$, where $x_\alpha=x'_0(\alpha)$ and $y_\alpha=x''_m(\alpha)$. 
Let $R_i(\delta)$ be the rectangle whose vertices are the points $x_i,\,x'_i(\delta),\,x''_i(\delta)$ and $x_{i+1}$ and let $B_i(\delta)$ be
the circular sector centered at $x_i$ and whose vertices are $x''_{i-1}(\delta)$ and $x'_i(\delta)$.

If we compute $\int_0^\delta L_{\xi_n}(\sigma_{\alpha})\,d\alpha$ it is not difficult to see that we get 
$$\int_0^\delta L_{\xi_n}(\sigma_{\alpha})\,d\alpha=
\sum_{i=0}^{m-1}\left(\int_{R_i(\delta)}\xi_nd\lcal^2\right)+\sum_{i=1}^{m-1}\left(\int_{B_i(\delta)}\xi_nd\lcal^2\right).$$
Moreover it holds $c_{\xi_n}(x^\alpha,y^\alpha)\leq L_{\xi_n}(\sigma^\alpha),$
hence we get
$$\int_0^\delta c_{\xi_n}(x^\alpha,y^\alpha)\,d\alpha\leq \sum_{i=0}^{m-1}\left(\int_{R_i(\delta)}\xi_nd\lcal^2\right)+
\sum_{i=1}^{m-1}\left(\int_{B_i(\delta)}\xi_nd\lcal^2\right).$$
If we pass to the limit as $n\to\infty$ we get, by using the uniform convergence of $c_{\xi_n}$ to $c$ on the left hand side 
and the $L^{q^*}$ convergence of $\xi_n$ to $\xi$ on the right hand  side,
$$\int_0^\delta c(x^\alpha,y^\alpha)\,d\alpha\leq \sum_{i=0}^{m-1}\left(\int_{R_i(\delta)}\xi d\lcal^2\right)+
\sum_{i=1}^{m-1}\left(\int_{B_i(\delta)}\xi d\lcal^2\right).$$
Then we divide by $\delta$ and we pass to the limit as $\delta\to 0$. Using the fact that $c$ is continuous we have
$$\lim_{\delta\to 0}\,\frac{1}{\delta}\int_0^\delta c(x^\alpha,y^\alpha)\,d\alpha=c(x,y).$$
On the other side, we may notice that the areas of the sectors $B_i$ may be estimated by $C\delta^2$ 
and hence we have, for $\delta\to 0$,
$$\frac{1}{\delta}
\sum_{i=1}^{m-1}\left(\int_{B_i}\xi d\lcal^2\right)\leq mC||\xi||\delta\to 0.$$
On the contrary the integrals over $R_i$, when divided by $\delta$, 
converge on the integrals on the segments $S_i$, which give exactly the integral over the curve $\sigma$, i.e. $L_\xi(\sigma)$. 
We have consequently
$$\lim_{\delta\to 0}\frac{1}{\delta}\left(\sum_{i=0}^{m-1}\left(\int_{R_i}\xi d\lcal^2\right)+
\sum_{i=1}^{m-1}\left(\int_{B_i}\xi d\lcal^2\right)\right)=L_{\xi}(\sigma)<c_{\xi}(x,y)+{1\over k}.$$
This gives
$$c(x,y)<c_{\xi}(x,y)+{1\over k}$$
and, $k$ being arbitrary, we also get $c\leq c_{\xi}$ and the thesis.
\end{proof}

\begin{lem}\label{existapprox}
Let us assume that $q<2$ and let $\xi$ be non-negative function belonging to $L^{q^*}$, then there exists a sequence $(\xi_n)_n\in C^0(\Omega),\,\xi_n\geq 0,\,
\xi_n\to\xi\,\mbox{ in }L^{q^*}$, such that $c_{\xi_n}$ converges uniformly to $\overline{c}_\xi$ on $\Omega$. 
\end{lem}

\begin{proof}
It is easy to see that for every $(x,y)\in\Omega^2$ there exists a sequence of non-negative continuous functions $(\xi_n)_n$ converging to $\xi$ in $L^{q^*}(\Omega)$ such that $c_{\xi_n}$ converges in $C^0$ and $\cbxi(x,y)=\lim_n c_{\xi_n}(x,y)$.  Let $I$ be a finite set,  $(x_i,y_i)\in \Omega^2$ for all $i\in I$ and for every $i$, let $(\xi_n^i)$ be  a sequence of non-negative continuous functions  converging to $\xi$ in $L^{q^*}(\Omega)$ such that $\cbxi(x_i,y_i)=\lim_n c_{\xi_n^i}(x_i,y_i)$. Let us set $\xi_n:=\max_{i\in I} \xi_n^i$, we then have $(\xi_n)_n$ converging to $\xi$ in $L^{q^*}(\Omega)$,  and 
\[\cbxi(x_i,y_i) \leq \liminf_n \; c_{\xi_n}(x_i,y_i) 
\leq \limsup_n \; c_{\xi_n}(x_i,y_i) \leq \cbxi(x_i,y_i).\]
We thus have $\cbxi(x_i,y_i)=\lim_n c_{\xi_n}(x_i,y_i)$ for every $i\in I$. Now,  let $(x_i,y_i)_{i\in \N}$ be a dense sequence of points of $\Omega^2$. From what preceeds, for every $n$, there exists a continuous non-negative $\xi_n$ such that 
\[\Vert \xi_n -\xi\Vert_{L^{q^*}}\leq \frac{1}{n}, \; \vert \cbxi(x_k,y_k)-c_{\xi_n}(x_k,y_k)\vert \leq \frac{1}{n},\; \forall k\leq n.\]
By the H\"older estimate of proposition \ref{cxicomp} and Ascoli's theorem, passing to a subsequence if necessary, we may assume that $c_{\xi_n}$ converges in $C^0$ to some $c$. Since obviously $c(x_k,y_k)=\cbxi(x_k,y_k)$ for all $k$ , we deduce $c=\cbxi$ and the desired result follows.  
\end{proof}

The next lemma enables us to extend $L_{\xi}$ in some sense when $\xi\geq 0$ is only $L^{q^*}$:

\begin{lem}\label{prolLxi}
Let us assume that $q<2$. Let $Q\in \QQq$, $\xi$ be a non-negative  element of $L^{q^*}$, and $(\xi_n)_n$ be a sequence of non-negative continuous functions that converges to $\xi$ in $L^{q^*}$, then we have the following:
\begin{enumerate}
\item[(i)] $(L_{\xi_n})_n$ converges strongly in $L^1(C,Q)$ to some limit which is independent of the  approximating sequence $(\xi_n)_n$ and which will again be denoted $L_\xi$. 
\item[(ii)] The following equality holds:
\begin{equation}\label{eglc}
\int_{\Omega} \xi(x) i_Q(x)\ dx=\int_{C} L_{\xi}(\sigma)\ dQ(\sigma).
\end{equation}
\item[(iii)] The following inequality holds for $Q$-a.e. $\sigma\in C$:
\begin{equation}\label{ineglc}
L_{\xi}(\sigma)\geq \oc_{\xi}(\sigma(0),\sigma(1)).
\end{equation}
\end{enumerate}
 
\end{lem}

\begin{proof}

For all $n$ and $m$ in $\N$ we have: 
\begin{eqnarray*}
\int_C \vert L_{\xi_n}(\sigma)-L_{\xi_m}(\sigma)\vert\ dQ(\sigma)\!&=& \!
\int_C \left\vert \int_0^1 \!\left(\xi_n(\sigma(t))\!-\!\xi_m(\sigma(t))\right)\vert\dot\sigma(t)\vert\ dt\right\vert\ \!dQ(\sigma)\\
&\le&\int_{\Omega} \vert \xi_n(x)-\xi_m(x)\vert i_Q(x)\ dx\\
&\le &  \Vert \xi_n-\xi_m\Vert_{L^{q^*}} \Vert i_Q\Vert_{L^q}.   
\end{eqnarray*}
This implies that $(L_{\xi_n})_n$ is a Cauchy sequence in $L^1(C,Q)$ and it is obvious, from the previous inequality, that its $L^1(C,Q)$ limit does not depend on the  approximating sequence $(\xi_n)_n$.

\smallskip
The proof of $(ii)$ follows from $(i)$:
\begin{eqnarray*}
\int_{\Omega}\xi(x)i_Q(x)\ dx&=&\lim_{n} \int_{\Omega}\xi_n(x)i_Q(x)\ dx\\
&=& \lim_{n}\int_C L_{\xi_n}(\sigma)\ dQ(\sigma)\\
&=&\int_C L_{\xi}(\sigma)\ dQ(\sigma).
\end{eqnarray*}

To prove $(iii)$ we choose an approximating sequence $(\xi_n)_n$ as in lemma \ref{existapprox} and pass to the limit in
\begin{equation*}
L_{\xi_n}(\sigma)\geq c_{\xi_n}(\sigma(0),\sigma(1)).\qedhere
\end{equation*}
\end{proof}

\subsection{Characterization of optimal transport strategies}\label{ots}

In this paragraph, our aim is to make the formal arguments of paragraph \ref{foc} rigorous in order to characterize optimal transport strategies. This can be done under the additional assumption that $H$ is strictly convex. First, we relate the optimality condition \pref{varineq} to the Monge-Kantorovich problem with cost $\oc_{\oxi}$:

\begin{prop}\label{varineqmk}
Let us assume that $q<2$ and that $H$ is strictly convex. If $\oQ$ solves \pref{lepbme} and $\oxi:=H'(i_{\oQ})$ then we have:
\begin{equation}\label{infegal}
\int_{\Omega}\oxi i_{\oQ}=\inf_{Q\in\QQq} \int_{\Omega}\oxi i_Q=\inf_{\gamma\in \Pi(\mu_0,\mu_1)} \int_{\omb\times \omb} \oc_{\oxi}(x,y)d\gamma(x,y).
\end{equation}

\end{prop}

\begin{proof}
Let us recall that from proposition \ref{foc}, we have:
\begin{equation}\label{rappfoc}
\int_{\Omega}\oxi i_{\oQ}=\inf_{Q\in\QQq} \int_{\Omega}\oxi i_Q.
\end{equation}

Let $\xi$ be a non-negative element of $L^{q^*}$ and let $Q\in\QQq$, using Lemma \ref{prolLxi} and the definition of $\QQq$ yields:
\[\begin{split}
\int_{\Omega} \xi i_Q  &=\int_C L_{\xi} dQ\geq \int_{C} \oc_{\xi}(\sigma(0)),\sigma(1)) dQ(\sigma) \\
& \geq \inf_{\gamma\in \Pi(\mu_0,\mu_1)} \int_{\omb\times \omb} \oc_{\xi}(x,y)d\gamma(x,y).
\end{split}\]
We then have, for all $\xi\in L^{q^*}$, $\xi\geq 0$:
\begin{equation}\label{ineglq}
\inf_{Q\in\QQq} \int_{\Omega}\xi i_Q \geq  \inf_{\gamma\in \Pi(\mu_0,\mu_1)} \int_{\omb\times \omb} \oc_{\xi}(x,y)d\gamma(x,y)
\end{equation}
and by a similar argument, for all $\xi \in C^0(\omb,\R_+)$:
\begin{equation}\label{inegc0}
\inf_{Q\in\QQ} \int_{\omb}\xi di_Q \geq  \inf_{\gamma\in \Pi(\mu_0,\mu_1)} \int_{\omb\times \omb} \oc_{\xi}(x,y)d\gamma(x,y).
\end{equation}

Let  $\xi \in C^0(\omb,\R_+)$, and $\e>0$, for every $x$ and $y$ in $\omb$, there exists $\sigma^{x,y}_\e\in C^{x,y}$ such that $(x,y)\mapsto \sigma^{x,y}_\e$ is measurable (see for instance \cite{cv}) and by Lemma \ref {ccoincid}
\begin{equation}\label{geodeps}
L_{\xi}(\sigma^{x,y}_\e)\leq c_{\xi}(x,y)+\eps=\oc_{\xi}(x,y)+\eps.
\end{equation}
Let $\gamma\in \Pi(\mu_0,\mu_1)$ and let us define the element of $\QQ$: $Q_\e:=\delta_{\sigma^{x,y}_\e}\otimes \gamma$, we  then have:
\[\int_{\omb} \xi di_{Q_\e} =\int_{\omb\times\omb} L_{\xi}(\sigma^{x,y}_\e)d\gamma(x,y)
\leq \int_{\omb\times\omb} \oc_{\xi}(x,y)d\gamma(x,y)+\eps.\]
Since $\gamma$ and $\eps$ are abitrary, using \pref{inegc0} we obtain
\begin{equation}\label{egc0}
\inf_{Q\in\QQ} \int_{\omb}\xi di_Q = \inf_{\gamma\in \Pi(\mu_0,\mu_1)} \int_{\omb\times \omb} \oc_{\xi}(x,y)d\gamma(x,y),\; \forall \xi\in C^0(\omb,\R_+).
\end{equation}

In what follows, for every $\mu\in \M_+ (\omb)$, we extend $\mu$ by $0$ outside $\omb$. Let $(\rho_n)_n$ be a standard mollifying sequence. For $n\in \N^*$ let us consider the regularized problem:
\begin{equation}\label{pbmen}
\inf_{Q\in \QQ} \int_{\R^2} H(\rho_n\star i_Q)
\end{equation}
The existence of a solution $\oQ_n$ of \pref{pbmen} can be obtained by similar arguments as in theorem \ref{existence} (using lemma \ref{tight} and the fact that the $L^1$ norm of $\rho_n\star i_Q$ equals the total mass of $i_Q$). Proceeding as in proposition \ref{foc} and defining $j_n:=\rho_n\star i_{\oQ_n}$, $\xi_n:=H'(j_n)$, $\eta_n:=\rho_n\star \xi_n$, we have:
\begin{equation}\label{varineqn}
\int_{\R^2} H'(\rho_n \star i_{\oQ_n})(\rho_n \star i_{\oQ_n})=\int_{\R^2} \xi_n j_n=\int_{\R^2} \eta_n d i_{\oQ_n}=\inf_{Q\in \QQ} \int_{\R^2} \eta_n d i_{Q}
\end{equation}
With \pref{egc0}, we then get:
\begin{equation}\label{mkn}
 \int_{\R^2} \eta_n d i_{\oQ_n}=\int_{\R^2} \xi_n j_n= \inf_{\gamma\in \Pi(\mu_0,\mu_1)} \int_{\omb\times \omb} \oc_{\eta_n}(x,y)d\gamma(x,y).
\end{equation}
By convexity of $H$, we also have:
\begin{equation}\label{inegenn}
\int_{\R^2} H(j_n)\leq \int_{\R^2} H(\rho_n\star i_{\oQ})\leq \int_{\R^2} \rho_n \star H(i_{\oQ})
\end{equation}
which implies that $j_n$ is bounded in $L^q$. Passing to subsequences, we may therefore assume:
\begin{equation}\label{cvgce0}
j_n \weakto j \mbox{ in $L^q$, } \xi_n \weakto \xi \mbox{ in $L^{q^*}$,} \; \eta_n \weakto \xi \mbox{ in $L^{q^*}$.}
\end{equation}
Since the total mass of $i_{\oQ_n}$ is the same of $j_n$ and $j_n$ is bounded in $L^q$ (and hence in $L^1$), we get a bound on $i_{\oQ_n}(\Omega)$ and, from lemma \ref{tight}, we may also assume:
\begin{equation}\label{cvgce1}
\oQ_n \ws Q \mbox{ in $\M_+(C)$, } i_{\oQ_n} \ws i \mbox{ in $\M_+(\omb)$.}
\end{equation}
It is obvious that $i=j$ and lemma \ref{I-IQ} implies $j\geq i_Q$. With \pref{inegenn} and the monotonicity of $H$, we then get:
\begin{equation}
\int_{\Omega} H(i_Q)\leq \int_{\Omega} H(j)\leq \liminf_n \int_{\R^2} H(j_n) \leq \int_{\Omega} H(i_{\oQ}).
\end{equation}
With the strict convexity of $H$ and the optimality of $\oQ$,  this also yields
\begin{equation}\label{limenergie}
i_{\oQ}=i_Q=j\in L^q \mbox{ and } \liminf_n \int_{\R^2} H(j_n)=\int_{\R^2}  H(i_{\oQ}).
\end{equation}
Up to some subsequence (as $\oxi=H'(i_{\oQ})\in L^{q^*}$), this also implies 
\[ H(j_n)-H(i_{\oQ})-\oxi (j_n-i_{\oQ})\to 0 \mbox{ a.e. and  in $L^1$}\]
and using the strict convexity of $H$, we deduce that $j_n$ converges a.e. to $i_{\oQ}$. This implies that $\xi_n$ converges a.e. to $\oxi=H'(i_{\oQ})$ and that $\xi=\oxi$. With Fatou's Lemma and \pref{mkn}, we therefore obtain:
\[\begin{split}
\int_{\Omega} \oxi i_{\oQ}& =  \int_{\R^2} H'(i_{\oQ}) i_{\oQ}\\
& \leq \liminf_n \int_{\R^2} \xi_n j_n\\
&= \liminf_n  \inf_{\gamma\in \Pi(\mu_0,\mu_1)} \int_{\omb\times \omb} \oc_{\eta_n}(x,y)d\gamma(x,y).
\end{split}\] 
Using $\eta_n \weakto \oxi$ and remark \ref{cxifaible}, from the uniform convergence of $c_{\eta_n}$ to a cost $c\leq c_{\oxi}$, we get
\[\int_{\Omega} \oxi i_{\oQ} \leq  \inf_{\gamma\in \Pi(\mu_0,\mu_1)} \int_{\omb\times \omb} \oc_{\oxi}(x,y)d\gamma(x,y)\]
together with \pref{rappfoc} and \pref{ineglq}, this completes the proof.
\end{proof}

The characterization  of  optimal transport strategies then reads as:

\begin{thm}\label{caractopt}
Let us assume that $q<2$ and that $H$ is strictly convex. A transportation strategy $(\ogm,\ovp)$ is optimal if and only if, setting $\oQ:=Q_{\ogm, \ovp}$ and $\oxi:=H'(i_{\oQ})$, one has:

\begin{enumerate}

\item $\ogm$ solves the Monge-Kantorovich problem:
\begin{equation}\label{foc1}
\inf_{\gamma\in \Pi(\mu_0,\mu_1)} \int_{\omb\times \omb} \oc_{\oxi}(x,y)d\gamma(x,y),\end{equation}

\item for $\oQ$-a.e. $\sigma\in C$, one has:
\begin{equation}\label{foc2}
L_{\oxi}(\sigma)=\oc_{\oxi}(\sigma(0),\sigma(1)).\end{equation}

\end{enumerate}

\end{thm}

\begin{proof}
Let us assume first that the transportation strategy $(\ogm,\ovp)$ is optimal and set $\oQ:=Q_{\ogm, \ovp}$ and $\oxi:=H'(i_{\oQ})$. From Proposition \ref{varineqmk} and Lemma \ref{prolLxi}, we get:
\[\begin{split}
 \int_{\omb\times \omb} \oc_{\oxi}(x,y)d\ogm(x,y) &=\int_C \oc_{\oxi}(\sigma(0),\sigma(1))d\oQ(\sigma)\\
  &\leq \int_C L_{\oxi} d\oQ=\int_{\Omega} \oxi i_{\oQ}\\
&= \inf_{\gamma\in \Pi(\mu_0,\mu_1)} \int_{\omb\times \omb} \oc_{\oxi}(x,y)d\gamma(x,y)
\end{split}\]
this proves that $\ogm$ solves \pref{foc1} and implies that the inequalities above are equalities. We therefore deduce \pref{foc2} from the inequality $ \oc_{\oxi}(\sigma(0),\sigma(1))\leq L_{\oxi}(\sigma)$ .

\smallskip

Conversely, assume that the transportation strategy $(\ogm,\ovp)$ satisfies the two conditions of the theorem. Condition \pref{foc2} firstly yields:
\[\int_{\Omega} \oxi i_{\oQ}=\int_C L_{\oxi} d\oQ=\int_C  \oc_{\oxi}(\sigma(0),\sigma(1))  d\oQ(\sigma)=
\int_{\omb\times \omb} \oc_{\oxi}(x,y)d\ogm(x,y)\]
Secondly, if $Q=Q_{\gamma,p}\in\QQq$, one has:
\[\int_{\Omega} \oxi i_{Q}=\int_C L_{\oxi} dQ \geq \int_C  \oc_{\oxi}(\sigma(0),\sigma(1))  dQ(\sigma)=
\int_{\omb\times \omb} \oc_{\oxi}(x,y) d\gamma(x,y)\]
and since $\ogm$ solves \pref{foc1}, we finally have
\[\int_{\Omega} \oxi i_{\oQ}\leq \int_{\Omega} \oxi i_{Q}, \; \forall Q\in \QQq\]
which, with proposition \ref{foc}, proves that $(\ogm,\ovp)$ is optimal.

\end{proof}

\begin{rem}
Let us remark that $\oxi=H'(i_{\oQ})$ (with $\oQ$ solving \pref{lepbme}) solves the following (dual  of \pref{lepbme}) problem:
\begin{equation}\label{pbdual}
\sup_{\xi \in L^{q^*}, \; \xi \geq 0} W(\xi)-\int_{\Omega} H^*(\xi(x))dx
\end{equation}
where $H^*$ is the Fenchel transform of $H$ and:
\[W(\xi):=\inf_{\gamma\in \Pi(\mu_0,\mu_1)} \int_{\omb\times \omb} \oc_{\xi}(x,y)d\gamma(x,y).\]
Indeed, it follows from proposition \ref{varineqmk} and $\oxi=H'(i_{\oQ})$ that
\[W(\oxi)-\int_{\Omega} H^*(\oxi(x))dx=\int_{\Omega} (\oxi i_{\oQ}-H^*(\oxi))=\int_{\Omega} H(i_{\oQ}).\]
If $\xi\in L^{q^*}$ is non-negative, we deduce from \pref{ineglq} and Young's inequality:
\[W(\xi)\leq \int_{\Omega} \xi i_{\oQ} \leq \int_{\Omega} H(i_{\oQ})+ \int_{\Omega} H^*(\xi)\]
which proves that $\oxi$ solves \pref{pbdual}.
\end{rem}

Let us see, through an easy exemple, an application of Theorem \ref{caractopt}.

\begin{ex}
Suppose that $\omb$ contains the two segments $A=\{0\}\times [0,1]$ and $B=\{1\}\times [0,1]$ and the square $S=[0,1]\times[0,1]$ which is their convex hull. Set $\mu_1=\huno\res A$ and $\mu_2=\huno\res B$ and denote by $T:A\to C$ the map that associates to every point $(0,x)\in A$ the curve $T(x)$ given by $T(x)(t)=(t,x)$, i.e. the horizontal segment from $A$ to $B$ starting from $x$. Set $Q=T\sharp\mu_1$. It is clear that $Q$ comes from an admissible transportation strategy linking $\mu_1$ to $\mu_2$ and it is not difficult to see that the traffic intensity $i_Q$ has constant density $1$ on $S$ and $0$ elsewhere. We consider two particular cases only: we claim that $Q$ is optimal if $\Omega=]0,1[\times]0,1[$ while it is not if $S$ is compactly contained in $\Omega$. Indeed, if $\omb=S$, the metric induced by $i_Q$ is the euclidean metric, the paths $T(x)$ are geodesic and the transport plan induced by $Q$ is optimal according to this metric. On the other hand, if $\Omega$ is larger than $S$, then all the segments $T(x)$ that are very close to the upper or lower boundary of $S$ are not geodesic accoding to this metric, because they could be improved by non-straight line paths which arrive up to zone $\Omega\setminus S$ where $i_Q=0$ and the trasportation is cheaper. In the former case, consequently, the sufficient optimality conditions are satisfied, while in the latter the geodesic conditions on the paths (Wardrop condition, see the nect section) is not and prevents optimality.
\end{ex}

\section{Application to equilibria of Wardrop type}\label{we}

In this final section, we relate the results of the previous sections to some concepts of equilibria of Wardrop type.  Modelling congestion as in paragraph \ref{tcm} enables us to extend the concept of Wardrop equilibrium to a continuous setting. 

\smallskip

Let us consider a congestion function $g$ : $\R_+\rightarrow \R_+$ which is continuous increasing and satisfies $az^{q-1}\leq g\leq b(z^{q-1}+1)$ for all $z\in \R_+$ and some $q\in(1,2)$ and non-negative constants $a$ and $b$. Then for any transportation strategy $(\gamma,p)$ such that $I_{\gamma,p}$ (defined by \pref{deftrafintgp}) belongs to $L^q$, the transportation cost function resulting from the strategy $(\gamma,p)$ is $\oc_{\xi}$ for $\xi:=g\circ I_{\gamma,p} \in L^{q^*}$. Roughly speaking, an equilibrium is then a  transportation  strategy $(\gamma, p)$ that satisfies Wardrop stability condition (i.e. $Q_{\gamma,p}$ gives full mass to the set of "geodesics" for the metric $\xi=g\circ I_{\gamma,p}$) and the additional requirement that $\gamma$ is an optimal transportation plan between $\mu_0$ and $\mu_1$ for the cost resulting from $(\gamma,p)$.  This leads to the following

\begin{defi}\label{defward}
A transportation strategy $(\ogm, \ovp)$ is said to be an equilibrium if $I_{\ogm,\ovp}\in L^q$ and, setting $\oxi:=g\circ I_{\ogm,\ovp}$ one has
\begin{enumerate}
\item $L_{\oxi}(\sigma)=\oc_{\oxi}(\sigma(0),\sigma(1))$ for $Q_{\ogm,\ovp}$-a.e. $\sigma\in C$, 

\item $\ogm$ solves the Monge-Kantorovich problem:
\[  \inf_{\gamma\in \Pi(\mu_0,\mu_1)} \int_{\omb\times \omb} \oc_{\oxi}(x,y)d\gamma(x,y).  \]
\end{enumerate}

\end{defi}

Only the first condition above is linked to Wardrop's original equilibrium concept. Imagine that some social planner chooses the transportation plan $\gamma$, then the second equilibrium condition expresses that $\gamma$ is optimal for the transportation cost resulting from $\gamma$ itself and the traveler's individual behavior. Our notion of equilibrium can therefore be viewed as a refinement of the Wardrop equilibrium or its generalization to the case where the transportation plan is not given a priori. 

\smallskip

A direct application of theorems \ref{existence} and \ref{caractopt} then gives the existence of equilibria together with a variational characterization:

\begin{thm}\label{cwe}
Under the assumptions of this paragraph, there exists an equilibrium. Moreover $(\ogm,\ovp)$ is an equilibrium if and only if $\oQ:=Q_{\ogm,\ovp}$ solves the minimization problem:
\begin{equation}
\inf_{Q\in \QQq} \int_{\Omega} H_g(i_Q(x))dx \quad\mbox{ with } H_g(z):=\int_0^z g(s)ds, \; \forall  z\in\R_+.
\end{equation}

\end{thm}

\begin{rem}
A slightly different situation, which can be relevant in some applications,  occurs when  the transportation plan $\ogm\in \Pi(\mu_0,\mu_1)$ is fixed and not only the marginals $\mu_0$ and $\mu_1$.  In this case,    one defines equilibria as the set of $\ovp$'s such that $(\ogm,\ovp)$ satisfies the first condition (Wardrop) of definition \ref{defward}. If the set:
\[\QQqg:=\{Q\in \M_+^1(C) \mbox{ : } (e_0,e_1)\sharp Q=\ogm,\; i_Q\in L^q\} \]
is nonempty, then slightly adapting our arguments, we have existence of equilibria and $\ovp$ is an equilibrium if and only if $\oQ:=Q_{\ogm,\ovp}$ solves the minimization problem:
\begin{equation*}
\inf_{Q\in \QQqg} \int_{\Omega} H_g(i_Q(x))dx.
\end{equation*}

\end{rem}



\end{document}